\documentclass[a4paper]{article}
\usepackage {amssymb,latexsym}
\usepackage {amsmath}
\newtheorem{thm}{Theorem}

\newtheorem{cor}{Corollary}

\newtheorem{de}{Definition}

\newenvironment{proof}{
                        \noindent{\bf\small Proof: }\small}
                                       {\hfill {$\mathbf \Box$}\medskip}

\newcommand{\K}{\mathbb{K}}
\newcommand{\N}{\mathbb{N}}
\newcommand{\Z}{\mathbb{Z}}

\newcommand{\nat}{\mathbb N}

\newcommand{\ra}{\rightarrow}

\newcommand{\korps}{\mathbb K}

\newcommand{\Ug}{\mathcal{U}\mathfrak{g}}

\makeatletter
\@addtoreset{equation}{section}

\makeatother

\title{On Graded Deformations of The universal Enveloping Algebra of a Color Lie Algebra}
\author{
\\Toukaiddine Petit\footnote{{\tt toukaiddine.petit@icloud.com}. }
}

\date{}

\begin{document}

\maketitle

\begin{abstract}
Let $\mathfrak{g}$ be a Color Lie Algebra and $\mathcal{U}(\mathfrak{g})$ its the universal Enveloping Algebra. We define the notion of graded deformations and we give explicit graded deformations of the universal Enveloping Algebra of $\mathfrak{g}$.
\end{abstract}
\section{Premilinaries}
Throughout this paper  groups are assumed to be abelian and $\korps$ is
a field of characteristic zero. We recall some notation for graded
algebras and graded modules \cite{NV}, and some facts on color Lie
algebras from  \cite{S1,S2,S3,SZ}.

\subsection{Color Lie algebras}
The concept of color Lie algebras is related to an abelian group $G$
and an anti-symmetric bicharacter $\varepsilon:G \times G
\rightarrow\K^\times$, i.e.,
 \begin{align}
&\varepsilon\left(g,h\right)\varepsilon\left(h,g\right)=1,\\
    &\varepsilon
\left(g,hk\right)=\varepsilon\left(g,h\right)\varepsilon\left(g,k\right),\\
    &\varepsilon\left(gh,k\right)=\varepsilon\left(g,k\right)\varepsilon\left(h,k\right),
\end{align}
where $g, h, k \in G$ and $\K^\times $ is the multiplicative group
of the units in $\K$.\par \vskip 5pt

 A $G$-graded space $L=\oplus_{g\in G} L_g$ is
said to be a $G$-graded $\varepsilon$-Lie algebra (or simply, color
Lie algebra), if it is endowed with a bilinear bracket
$\left[-,-\right]$ satisfying the following conditions
\begin{equation}
\left[ L_g, L_h \right]\subseteq L_{gh},
\end{equation}
\begin{equation}
    \left[a,b\right]=-\varepsilon\left(|a|,|b|\right)\left[b,a\right],
\end{equation}
\begin{equation}
\varepsilon\left(|c|,
|a|\right)\left[a,\left[b,c\right]\right]+\varepsilon\left(|a|,|b|\right)\left[b,\left[c,a\right]\right]+
\varepsilon\left(|b|, |c|\right)\left[c,\left[a,b\right]\right]=0,
\end{equation}
where $g, h \in G$,  and $a, b, c \in L$ are homogeneous
elements.\par \vskip 5pt

 For example, a
super Lie algebra is exactly  a $\Z_2$-graded $\varepsilon$-Lie
algebra where
\begin{equation}
    \varepsilon(i,j)=(-1)^{ij},\quad \forall\quad i,j\in\Z_2.
\end{equation}
\subsection{Graded Cohomology of Cartan Eilenberg of Color Lie Algebra}
A vector space $M$ is $G$-graded if the is a family of subspaces $(M_g){g\in G}$ such that $M=\oplus M_g$ . Let $M,N$ two graded vector spaces. A linear map from $M$ to $N$ is of homogeneous degree $\alpha$ if $f(M_{\beta})\subset N_{\alpha+\beta}$ for $\beta$.  Denote by $Hom_{\alpha}(M,N)$ the set the linear maps of homogeneous degree $\alpha$. The graded vector space $Hom_{gr}(M,N)=\oplus Hom_{\alpha}(M,N)$ equipped by the bracket defined by $[f,g]=f\circ g-\varepsilon(f,g)g\circ f$ is a color Lie algebra. Let $L$ be a color Lie algebra. A  graded representation of $L$ in $M$ is a linear map of degree zero, $\rho: L\rightarrow Hom_{gr}(M,M)$ satisfying $[\rho(f),\rho(g)]=\rho([f,g])$. One said that $M$ is a graded  L module. In general a $n$ linear map $f:L\times\cdots\times L\rightarrow M$ is of homogeneous degree $\alpha$ if $f(X_{\alpha_1},\cdots,X_{\alpha_n})$ is homogeneous of degree $\alpha+\alpha_1+\cdots +\alpha_n$.  Denote by $Hom^n_{\alpha}(L\times\cdots\times L, N)$ the set the linear maps of homogeneous degree $\alpha$. The $Hom^n_{gr}(L\times \cdots\times L,N)=\oplus Hom^n_{\alpha}(L\times\cdots\times L, N)$ is a $G$ graded vector space. It admits a $G$ graded L module given by $\rho(X)(f)(X_{\alpha_1},\cdots,X_{\alpha_n})=$
\begin{equation}
\rho(X)(f(X_{\alpha_1},\cdots,X_{\alpha_n}))-\sum^n_{i=1}\varepsilon(\alpha,\gamma+\alpha_1+\cdots \alpha_{i-1})f(X_{\alpha_1},\cdots, [X,X_{\alpha_r}], X_{\alpha_n})
\end{equation}
An element of $f$ of homogeneous of degree $\gamma$ is called $\varepsilon$  skew symmetric if

\begin{equation}
f(X_{\alpha_1},\cdots,X_{\alpha_i},\cdots, X_{\alpha_j}, X_{\alpha_n})=-\varepsilon(\alpha_i,\alpha_j)f(X_{\alpha_1},\cdots,X_{\alpha_j},\cdots, X_{\alpha_i}, X_{\alpha_n})
\end{equation}

for $i<j$ 
We set \\
$C^0(L,N)=N$,  $C^1_{gr}(L,N)= Hom_{gr}(L,N)$, $C_{gr}^n(L,N)$ the set of elements which are $\varepsilon$ skew symmetric of $Hom^n_{gr}(L\times \cdots\times L,N)$.  for $n>1$. The graded space called $n$ graded cochains.
The linear coboudary operator 
\begin{equation}
d^n:C_{gr}^n(L,N)\rightarrow C_{gr}^{n+1}(L,N)
\end{equation}
defined by for $n>1 $

\begin{align}
&\delta^n(f)\left(x_1, \cdots, x_{n+1}\right)\\
&=\sum_{i=1}^{n+1} (-1)^{i+1}\varepsilon(\gamma,\alpha_i)\varepsilon_i\;x_i\cdot f\left( x_1, \cdots, \hat{x_i}, \cdots, x_{n+1}\right)\\
&\quad +\sum_{1\leq i< j \leq n+1} (-1)^{i+j}\varepsilon(\gamma,\alpha_i)\varepsilon(\gamma,\alpha_j)\varepsilon_i
\varepsilon_j\;f\left( [x_i, x_j], x_1,
\cdots, \hat{x_i}, \cdots, \hat{x_j}, \cdots, x_{n+1}\right),
\end{align}

for all $f \in  C_{gr}^n(L,N)_{\gamma}$, with $\varepsilon_i=\prod_{h=1}^{i-1}\varepsilon(|x_h|, |x_i|)$ $i\geq 2$, $\varepsilon_1=1$ and the sign\quad$\widehat{}$\quad indicates that the element below it must be omitted. We have $\delta^{n+1}\circ \delta^n =0$, \cite{CP,S3}
\subsection{Graded Deformations of color Lie Algebras}
The field $K$ admits a $ G$ graded by setting $K_{e}=K$ and $K_\gamma=0$ if $\gamma\neq 0$. Then the ring of series with coefficient in $K$, $K[[t]]$ is $G$ graded by extension.  Let $E$  be a $G$-graded vector space. Then the $K[[t]]$ module $E[[t]]$ is $G$-graded module. Let $L$ be a Color Lie algebra, The $K[[t]]$ module $L[[t]]$ is a Color Lie algebra by extension.
\begin{de} Let $L$ be a Color Lie algebra over $K$ with multiplication $\phi_0$.
\begin{enumerate}
\item  A $G$ graded deformation of the Color Lie algebra $(L[[t]],\phi_0)$ (of degree zero) is a $G$ graded $K[[t]]$-module $L[[t]]$ equipped of a Color multiplication $\phi=\sum_{t}\phi_n t^n$ with 
$\phi_n\in C^2(L,L)_e $
\item  Two graded deformations $(L[[t]],\phi_1)$ and $(L[[t]],\phi_2)$ of $L[[t]],\phi_0)$ are said equivalent if they are isomorphic i.e. there is a graded automorphism $f=id+t f_1+t^2f_2+\cdots + t^n f^n+\cdots $ with $f\in C^1(L,L)_e$ such that $\phi_2=f^{-1}\circ \phi_1 \circ f\times f$.
\item A graded deformation $(L[[t]],\phi)$ is said trivial if all graded deformation is equivalent to $(L[[t]],\phi_0)$.
\item $L$ is called graded rigid if all graded deformation of $(L[[t]],\phi_0)$ is trivial.

\end{enumerate}
\end{de}
\begin{thm} Let $L$ be a color Lie algebra. If the second graded cohomology $H^2(L,L)_e$ is equal to zero, then $L$ is graded rigid.
\end{thm}
\subsection{Graded Hochschild cohomology}
Let $G$ be an abelian group with identity element $e$. We will
write $G$ as an multiplicative group.\par \vskip 5pt

An associative algebra $A$ with unit $1_A$, is said to be
$G$-graded, if there is a family $\left\{A_g|  g\in G\right\}$ of
subspaces  of $A$ such that $A=\oplus_{g \in G} A_g$ with  $1_A
\in A_e$ and $A_gA_h \subseteq A_{gh}$, for all $g, h \in G$. Any
element  $a \in A_g$ is  called homogeneous of degree $g$,
 and we write $|a|=g$. \par \vskip 5pt

A (left) graded $A$-module $M$ is a left $A$-module with an
decomposition $M=\oplus_{g \in G}M_g$ such that $A_g.M_h \subseteq
M_{gh}$. Let $M$ and $N$ be graded $A$-modules. Define
\begin{equation}
    \mathrm{Hom}_{A\mbox{-}{\rm gr}}(M,N)=
    \left\{f\in\mathrm{Hom}_{A}(M,N)|\ \  f(M_g)\subseteq N_g, \quad \forall\ \ g\in G\right\}.
\end{equation}

Let us recall the notion of graded Hochschild cohomology of a graded algebra $A$.
A graded $A$-bimodule is a $A$-bimodule $M=\oplus_{g \in G}M_g$ such that
 $A_g.M_h.A_k \subseteq M_{ghk}$. 
 The Hochschild graded cochain complex of A is  
$$ A \rightarrow  Hom_{gr}(A,A)\rightarrow  Hom_{gr}(A\times A,A)\rightarrow .. \rightarrow Hom_{gr}(A\times\cdots \times A,A)\rightarrow \cdots $$
whose diﬀerential is given by
\begin{equation}
d(f)(a_0,a_1,\cdots, a_p)=a_0 f(a_1,\cdots, a_p)-\sum_{i=0}^{p-1}(-1)^{i}f(a_0,a_1,\cdots, a_ia_{i+1},\cdots, a_p)+(-1)^pf(a_0,\cdots, a_{p-1})a_p
\end{equation}

\subsection{Enveloping Algebra of a Color Lie algebra}
Let $L$ be a color Lie algebra as above and $T(L)$ the tensor
algebra of the underlying $G$-graded vector space $L$. It is
well-known that $T\left(L\right)$ has a natural $\Z\times
G$-grading which is fixed by the condition that the degree of a
tensor $a_1\otimes...\otimes a_n$ with $a_i\in L_{g_i}$, $g_i\in
G$, $1\leq i\leq n$, is equal to $\left(n,g_1\cdots g_n\right)$.
The subspace of $T\left(L\right)$ spanned by homogeneous tensors
of order $\leq n$ will be denoted by $T^n\left(L\right)$. Let
$J\left(L\right)$ be the $G$-graded two-sided ideal of
$T\left(L\right)$ which is generated by
\begin{equation}
    a\otimes b-\varepsilon\left(|a|, |b|\right)b\otimes a-\left[a,b\right]
\end{equation}
with homogeneous $a, b \in L$. The quotient algebra
$U\left(L\right):=T\left(L\right)/J\left(L\right)$ is called the
universal enveloping algebra of the color Lie algebra $L$. The
$\K$-algebra $U\left(L\right)$ is a $G$-graded algebra and  has a
positive filtration by putting $U^n\left(L\right)$ equal to the
canonical image of $T^n\left(L\right)$ in $U\left(L\right)$.\par
\vskip 5pt
$I\left(L\right)$ be the $G$-graded two-sided ideal of
$T\left(L\right)$ which is generated by
\begin{equation}
    a\otimes b-\varepsilon\left(|a|, |b|\right)b\otimes a
    \end{equation}
with homogeneous $a, b \in L$. The quotient algebra
$S\left(L\right):=T\left(L\right)/I\left(L\right)$ is called the
universal symetric algebra of the color Lie algebra $L$. 
 
The canonical map $\verb"i":L\rightarrow U\left(L\right)$ is a
$G$-graded homomorphism and satisfies

\begin{equation}
\verb"i" \left(a\right)\verb"i" \left(b\right)-\varepsilon\left(|a|,
|b|\right)\verb"i"\left(b\right)\verb"i"\left(a\right)=\verb"i"\left(\left[a,b\right]\right).
\end{equation}

The $\Z$-graded algebra $G(L)$ associated with the filtered
algebra $U\left(L\right)$ is defined by letting
$G^n\left(L\right)$ be the vector space
$U^n\left(L\right)/U^{n-1}\left(L\right)$ and $G\left(L\right)$
the space $\oplus_{n\in\N}G^n\left(L\right)$ (note
$U^{-1}\left(L\right):=\left\{0\right\}$). Consequently,
$G\left(L\right)$ is a $\Z\times G$-graded algebra. The well-known
generalized Poincar\'{e}-Birkhoff-Witt theorem, , states
that the canonical homomorphism $i:L\rightarrow
U\left(L\right)$ is an injective $G$-graded homomorphism;
moreover, if $\left\{x_i\right\}_I$ is a homogeneous basis of $L$,
where the index set $I$ well-ordered. Set
$y_{k_j}:=\verb"i"\left(x_{k_j}\right)$, then the set of ordered
monomials $y_{k_1}\cdots y_{k_n}$ is a basis of $U\left(L\right)$,
where $k_j\leq k_{j+1}$ and $k_j< k_{j+1}$ if
$\varepsilon\left(g_j,g_j\right)\neq 1$ with $x_{k_j}\in L_{g_j}$
for all $1\leq j\leq n,n\in\N$. In case $L$ is finite-dimensional
$U\left(L\right)$ is a two-sided (graded) Noetherian algebra
.  The Poincar\'{e}-Birkhoff-Witt theorem shows that there is a canonical graded algebra isomorphism from $S(L)$ to $G(L)$.
\section{Graded Deformations of Enveloping algebra of A Color Lie Algebra}
\subsection{Graded Deformation by the Lie Color algebra}
\begin{thm} Let $\mathfrak{g}$ be a Color Lie Algebra and $\mathcal{U}(\mathfrak{g})$ its the universal Enveloping Algebra.
If  $(\mathfrak{g}[[t]],\mu_{t})$ is a non trivial graded deformation of $\mathfrak{g}$, then it is uniquely extended to a graded deformation of the associative algebra $\mathcal{U}(\mathfrak{g})$. 
\end{thm}

\begin{proof}
Let  be a 
$(\mathfrak{g}[[t]],\mu_{t})$ graded deformation of $\mathfrak{g}$ which is not  non trivial with $\mu_{t}=\sum_{n=0}^{\infty}\mu_{n}t^{n}$ such that the class of 
$\mu _{1}$ is not null in 
$\mathbf{H}^{2}(\mathfrak{g},\mathfrak{g})_e$.
Since $\mathfrak g$ is finite dimension, then the $\korps[[t]]$-module
$\mathfrak{g}[[t]]$ is graded 
isomorphic to the free module $\mathfrak{g}\otimes_{\korps}\korps[[t]]$.
Let  $e_1,\ldots,e_n$ be a homogeneous basis of
$\mathfrak{g}$, $y_1,\ldots,y_n$ the vector of images of the basis in 
$\mathcal{U}(\mathfrak{g})g$ and   $y'_1,\ldots,y'_n$ the vector of images of the basis in $\mathcal{U}\big(\mathfrak{g}[[t]]\big)$. Ones denotes $\bullet$
the multiplication of $\mathcal{U}\big(\mathfrak{g}[[t]]\big)$.
For all increase finite sequence $I=(i_1,\ldots,i_k)$ of indices in
$\{1,\ldots,n\}$
set $y_I:=y_{i_1}\cdots y_{i_k}$ in  $\Ug$ and 
$y'_I:=y'_{i_1}\bullet\cdots\bullet y'_{i_k}$ in
$\mathcal{U}\big(\mathfrak{g}[[t]]\big)$.
The  theorem of Poincar\'{e}-Birkhoff-Witt (is valid in the 
situation of Color Lie algebra which is a free module over a commutative ring, 
the $y_I$ forme a basis of  $\mathcal{U}(\mathfrak{g})$ over $\korps$ and 
$y'_I$ forme a basi of $\mathcal{U}\big(\mathfrak{g}[[t]]\big)$
over $\korps[[t]]$. We see the map $\Phi:Ug\otimes_\korps
\korps[[t]]\ra \mathcal{U}\big(\mathfrak{g}[[t]]\big)$ defined by 
$\Phi(y_I):=y'_I$ is graded isomorphic of $\korps[[t]]$-modules.
Let $\pi_t:\Ug\otimes_\korps \korps[[t]]~\times~Ug\otimes_\korps
\korps[[t]]~\ra \mathcal{U}(\mathfrak{g})\otimes_\korps \korps[[t]]$ the multiplication $\mathcal{U}(\mathfrak{g})\otimes_\korps\korps[[t]]$ induces by $\bullet$ et $\Phi$, i.e., $\pi_t(a,b):= \Phi^{-1}\big(\Phi(a)\bullet\Phi(b)\big)$. The restriction
of $\pi_t$ of elements of  $\mathcal{U}(\mathfrak{g})\times \mathcal{U}(\mathfrak{g})$ defined an application
$\korps$-bilinear  $\mathcal{U}(\mathfrak{g})\times \mathcal{U}(\mathfrak{g})\ra \mathcal{U}(\mathfrak{g})\otimes_\korps \korps[[t]]\subset
\mathcal{U}(\mathfrak{g})[[t]]$ that denote again by  $\pi_t$, i.e, $\pi_t(u,v)=\sum_{n=0}^{\infty}t^{n}\pi_{n}(u,v)$ all elements $u,v\in
\mathcal{U}(\mathfrak{g})$ o\`{u} $\pi_{n}\in\mathbf{Hom}_{gr}(\mathcal{U}(\mathfrak{g})\otimes \mathcal{U}(\mathfrak{g}),\mathcal{U}(\mathfrak{g}))$.
L'associativity of $\pi_t$ on three elements $u,v,w\in \mathcal{U}(\mathfrak{g})$ imply the equations $\sum_{s=0}^r\Big(\pi_s\big(\pi_{r-s}(u,v),w\big)-
\pi_s\big(u,\pi_{r-s}(v,w)\big)\Big)=0$ all $r\in\nat$.
Then, $\pi_t$ defines a bilinear multiplication
$\korps[[t]]$-associative on the $\korps[[t]]$-module
$\mathcal{U}(\mathfrak{g})\hspace{0.5mm}[[t]]$
(which is containing $\mathcal{U}(\mathfrak{g})\otimes_\korps \korps[[t]]$ as sub-module dense by
the topology $t$-adic) of usually manner 
\[
  \pi_t\left(\sum_{s}^\infty t^s u_s,\sum_{s'=0}^\infty t^{s'} v_{s'}\right)
      := \sum_{r=0}^\infty t^r \sum_{\stackrel{s,s',s''\geq 0}{s+s'+s''=r}}
                     \pi_{s''}(u_s,v_{s'})
\]
In particular, the map $\pi_0$ defines the associative multiplication
on the vector space $\mathcal{U}(\mathfrak{g})$, and
$(\mathcal{U}(\mathfrak{g})\hspace{0.5mm}[[t]],\pi_t)$ is a graded 
 associative deformation of $(\mathcal{U}(\mathfrak{g}),\pi_0)$. All finite increasing sequence  $I,J$ we have
$\pi_0(y_I,y_J)=\Phi^{-1}(y'_I\bullet y'_J)|_{t=0}$: by `well ordonning'
the product $y'_I\bullet y'_J$ we obtain of  linear combinaisons  of
$y'_K$ o\`{u} the increasing  finite sequence   $K$ is length less than or equal to the sum of the length of $I$ and $J$ which the coefficients do not contain only the structure constants of bracket $\mu_0$ of the Color Lie algebra $\mathfrak{g}$ since $t=0$. Then $\pi_0$ is equal the multiplication of the enveloping algebra $\mathcal{U}(\mathfrak{g})$ et $(\mathcal{U}(\mathfrak{g})[[t]],\pi_t)$ is a graded associative deformation of the enveloping algebra $\mathcal{U}(\mathfrak{g})$

It follows that $\pi_{1}$ is a graded $2$-cocycle of Hochschild of $\mathcal{U}(\mathfrak{g})$, and
the restriction of $\pi_{1}$ sur $X,Y\in\mathfrak{g}$ verifies
\begin{equation}
    \mu_{1}(X,Y)=\pi_{1}(X,Y)-\varepsilon(X,Y)\pi_{1}(Y,X) ~~~ \forall X,Y\in \mathfrak{g}.
\end{equation}
since the Color Lie algebra   $(\mathfrak{g}[[t]],\mu_t)$ is a Color Lie sub-algebra 
of $\mathcal{U}(\mathfrak{g}[[t]])$ which can be considered as a
associative subalgebra  of $(\mathcal{U}(\mathfrak{g})[[t]],\pi_t)$.
We assume that this deformation of $\mathcal{U}(\mathfrak{g})$ is trivial,  implies that there is
a formalle  graded isomorphism 
$\varphi_{t}=\sum_{r=0}^{\infty}\varphi_{r}t^{r}$, with
$\varphi_{0}=Id_{\mathcal{U}(\mathfrak{g})}$ et
$\varphi_{n}\in \mathbf{Hom}(\mathcal{U}(\mathfrak{g}),\mathcal{U}(\mathfrak{g}))_e$ such that
\[
    \varphi_{t}(\pi_{t}(u,v))=
            \pi_{t}(\varphi_{t}(u),\varphi_{t}(v)),
            ~~~\forall u,v\in \mathcal{U}(\mathfrak{g}),
\]
ce qui est \'{e}quivalent \`{a}
\begin{equation}
    \sum_{r=0}^{\infty}t^{r}\sum_{\stackrel{a,b\geq 0}{a+b=r}}
       \varphi_{a}(\pi_{b}(u,v))
         =\sum_{r=0}^{\infty }t^{r}\sum_{\stackrel{a,b,c\geq0}{a+b+c=n}}
            \pi_{a}(\varphi_{b}(u),\varphi_{c}(v))\quad \quad
                  \forall u,v\in \mathcal{U}(\mathfrak{g}).
\end{equation}
For $r=1$, this relation gives
\begin{equation}
   \pi_{1}(u,v)=
     (\delta_{H}\varphi_{1})(u,v)\quad \quad \forall u,v\in\mathcal{U}(\mathfrak{g})
\end{equation}
where $\delta_{H}$ is the differential operator  of Hochschild 
associated with the multiplication $\pi_0$ of the enveloping algebra.

Then the formulas  gives
\begin{equation}
   \mu_{1}(X,Y)  =  (\delta_{H}\varphi _{1})(X,Y)
                  -\varepsilon(X,Y)(\delta_{H}\varphi_{1})(Y,X) \nonumber 
                  \end{equation}
\begin{equation}
  =  X\varphi_1(Y)-\varphi(XY)+\varphi(X)Y
                      -\varepsilon(X,Y) (-Y\varphi_1(X)+\varphi(YX)-\varphi(Y)X) \nonumber              
\end{equation}

\begin{equation}
  =  (\delta _{CE}\varphi _{1})(X,Y)\quad \quad
                         \forall X,Y\in\mathfrak{g}
\end{equation}

where $\delta_{CE}$ is  the graded differential operator of Chevalley-Eilenberg

Then $\mu_{1}$ est is coboundary of Graded Chevalley-Eilenberg and its
class is null in $\mathbf{H}^{2}(\mathfrak{g},\mathfrak{g})_e$,
gives a contradiction.
\end{proof}
\begin{cor}Let $\mathfrak{g}$ be a Color Lie Algebra and $\mathcal{U}(\mathfrak{g})$ its the universal Enveloping Algebra. If the associative algebra $\mathcal{U}(\mathfrak{g})$ is graded rigid, then the Color Lie algebra $\mathfrak{g}$  is graded rigid.
\end{cor}
\subsection{Deformation by central extension }
Let $\mathfrak{g}$ be a Color Lie Algebra and $\omega$ a two graded cocycle with value in $K$ such that its class in non null in $H^2(\mathfrak{g},K)_e$. The central extension of $\mathfrak{g}$ by $c$ is defined by $\mathfrak{g}_{\omega}=\mathfrak{g}\oplus K$ and $[X+\alpha c,Y+\beta c]=[X,Y]+\omega(X,Y)c$
\begin{thm}  Let $\mathfrak{g}$ be a Color Lie Algebra and $\mathcal{U}_{\mathfrak{g}}$ its the universal Enveloping Algebra.
If the second graded cohomology with in $K$, $H_(\mathfrak{g},K)_e$ is different from zero, then the Enveloping algebra $\mathcal{U}_{\omega}\mathfrak{g}$ admits a non trivial graded deformation.
\end{thm}
\begin{proof}
Let $\omega\in\mathbf{Z}_{CE}^2(\mathfrak{g},\korps)_e$ be a graded $2$-cocycle of  non null class and $\mathfrak{g}_{t\omega}[[t]]$ the one dimension central extension of the Color lie algebra $\mathfrak{g}[[t]]=\mathfrak{g}\otimes_\korps \korps[[t]]$ over $K=\korps[[t]]$
 associated with $t\omega$ . 
In the enveloping algebra  $\mathcal{U}(\mathfrak{g}_{t\omega}[[t]])$ of
 $\mathfrak{g}_{t\omega}[[t]]$ we denote the multiplication by $\bullet$ and we
 consider the bilateral ideal 
 $\mathcal{I}:=(1-c')\bullet\mathcal{U}(\mathfrak{g}_{t\omega}[[t]])=
 \mathcal{U}(\mathfrak{g}_{t\omega}[[t]])\bullet (1-c')$ (where $c'$ is
 the image of $c$ in $\mathcal{U}\big(\mathfrak{g}_{t\omega}[[t]]\big)$)
 and the quotient algebra 
 $\mathcal{U}_{t\omega}\mathfrak{g}
 :=\mathcal{U}(\mathfrak{g}_{t\omega}[[t]])/\mathcal{I}$.
Let $e_1,\ldots,e_n$ a homogenous basis of
$\mathfrak{g}$ over $\korps$.  Then $c,e_1,\ldots,e_n$ is
a homogenous basis of $\mathfrak{g}_{t\omega}[[t]]$ sur $\korps[[t]]$.
Let  $y_1,\ldots,y_n$ the images of basis vecteurs  in
$\Ug$ et  $c',y'_1,\ldots,y'_n$ the images of basis  vecteurs in 
$\mathcal{U}\big(\mathfrak{g}_{t\omega}[[t]]\big)$.
For all finite increasing sequence  $I=(i_1,\ldots,i_k)$ of indices in 
$\{1,\ldots,n\}$
set $y_I:=y_{i_1}\cdots y_{i_k}$ in  $\Ug$ and 
$y'_I:=y'_{i_1}\bullet\cdots\bullet y'_{i_k}$ in 
$\mathcal{U}\big(\mathfrak{g}[[t]]\big)$.
The theorem of Poincar\'{e}-Birkhoff-Witt the $y_I$ form a basis of $\Ug$ over $\korps$, and the 
$c^{'\bullet i_0}\bullet y'_I$ (o\`{u} $i_0\in\nat$ et $c^{'\bullet i_0}:=1$)
form a basis of $\mathcal{U}\mathfrak{g}_{t\omega}[[t]]\big)$
over $\korps[[t]]$. In the quotient algebra 
$\mathcal{U}_{t\omega}\mathfrak{g}$, we see that
$c^{'\bullet i_0}$ can be identified to $1$. By denoting the multiplication in
$\mathcal{U}_{t\omega}\mathfrak{g}$ by $\cdot$ and the images of
$y'_1,\ldots,y'_n$ by the canonical projection by $y''_1,\ldots,y''_n$,
wee for all  $y'_I$ is sending to 
$y''_I:=y''_{i_1}\cdot\ldots\cdot y''_{i_n}$. It follows that all the 
elements $y''_I$ forme a basis of the quotient algebra 
$\mathcal{U}_{t\omega}\mathfrak{g}$.  As in the demonstration of the pervious theorem , the map  $\Phi:\Ug\otimes_\korps \korps[[t]]\ra
\mathcal{U}_{t\omega}\mathfrak{g}$ given by $y_I\mapsto y''_I$
defines  an isomorphism of free $\korps[[t]]$-modules.  In the similarly manner as the precedent demonstration we show  the multiplication induced over
$\Ug\otimes_\korps \korps[[t]]$ by the
multiplication $\cdot$ of $\mathcal{U}_{t\omega}\mathfrak{g}$ and 
$\Phi$ définies a suite of maps
$\pi_{t}=\sum_{r=0}^{\infty }\pi_{r}t^{r}$, o\`{u} $\pi_{r}\in
\mathbf{Hom} (\Ug\otimes\Ug,\Ug_e$ with the following property:
(1.) $\pi_t$  defines a graded associative deformation  $(\Ug,\pi_0)$,
i.e, a  $\korps[[t]]$-bilinear multiplication over
the graded  $\korps[[t]]$-module $\Ug\hspace{0.5mm}[[t]]$ (wich containing 
the $\korps[[t]]$-module $\Ug\otimes_\korps \korps[[t]]$ as sub-espace 
dense for the topology  $t$-adic), et (2.) $\pi_0$  is usually multiplication 
of the enveloping algebra  $\Ug$ of $\mathfrak{g}$.
It follows that  $\pi_{1}$ is graded  $2$-cocycle of Hochschild of $\Ug$, and over
all $X,Y\in \mathfrak{g}\subset \Ug$ we have the 
relation: $\omega(X,Y)1=\pi_{1}(X,Y)-\varepsilon(X,Y)\pi_{1}(Y,X)$ since
the Color Lie $\mathfrak{g}_{t\omega}[[t]]$ injecte in
the quotient algebra $\mathcal{U}_{t\omega}\mathfrak{g}$, and so  in
$\Ug\otimes_{\korps}\korps[[t]]\subset \Ug\hspace{0.5mm}[[t]]$.

We suppose that $\Ug$ is  rigid, then in particular, the deformation
$\pi_t$ is trivial. Then there exists a graded $1$-cocycle of Hochschild
$\varphi_{1}\in \mathbf{C}_{H}^{1}(\Ug,\Ug)$ such that
 $\pi_{1}=\delta_{H}(\varphi_{1})$. It follows that
$\forall X,Y\in \mathfrak{g}$:

\[
 \omega (X,Y)1=\pi_{1}(X,Y)-\varepsilon(X,Y)\pi_{1}(Y,X)
                  \]

\[=\delta_{H}(\varphi_{1})(X,Y)-\varepsilon(X,Y)
                   \delta_{H}(\varphi_{1})(Y,X)
                   =\delta_{CE}(\varphi_{1})(X,Y).
                   \]
Then $\omega$ is a graded cobord of Chevalley-Eilenberg and its
class is null in $\mathbf{H}_{CE}^{2}(\mathfrak{g},\korps)_e$,
gives a contradiction.
\end{proof}
\subsection{Deformation by Poincar\'{e}-Birkhoff-Witt Theorem }
\begin{de}Let $\chi$ be bi-character of $G$. Let $A=\oplus A_g$ a graded $\chi$ commutative associative algebra.  A Color Poisson bracket is a Color Lie bracket satisfying the Leibniz property $\lbrace a_{\alpha},a_{\beta}a_{\gamma}\rbrace=\lbrace a_{\alpha},a_{\beta}\rbrace a_{\gamma}+\chi(\alpha,\beta)a_{\beta}\lbrace a_{\alpha},a_{\gamma}\rbrace$. The algebra is called Color Poisson algebra.
\end{de}
The universal enveloping algebra $\mathcal{U}(\mathfrak{g})$ of the Color Lie algebra $\mathfrak{g}$ equipped by the canonical filtration $(\mathcal{U}_n(\mathfrak{g}))_{n}$  satisfying $\mathcal{U}(\mathfrak{g})_n.\mathcal{U}(\mathfrak{g})_m\subset \mathcal{U}(\mathfrak{g})_{n+m}$. The graded associative algebra $gr(\mathcal{U}(\mathfrak{g}))=\oplus gr^n(\mathcal{U}(\mathfrak{g}))$ for this filtation, is defined by $ gr^n(\mathcal{U}(\mathfrak{g}))=(\mathcal{U}(\mathfrak{g}))_n/(\mathcal{U}(\mathfrak{g}))_{n-1}$ with $(\mathcal{U}(\mathfrak{g}))_{-1}=0$. Let $\pi_n:\mathcal{U}(\mathfrak{g}))_n\rightarrow  gr^n(\mathcal{U}(\mathfrak{g}))$ the canonical projection. The multiplication of the algebra $gr(\mathcal{U}(\mathfrak{g}))$ is defined by
$$\pi_n(u)\pi_m(v):=\pi_{n+m}(u.v)$$
with $u\in \mathcal{U}(\mathfrak{g})_n$ and $v\in \mathcal{U}(\mathfrak{g})_m$
\begin{thm}Let $(\mathfrak{g},[,])$ be a Color Lie algebra, $\mathcal{U}(\mathfrak{g})$ its universal enveloping algebra. The associated graded algebra $gr(\mathcal{U}(\mathfrak{g}))$ is endowed of a Color Poisson bracket defined by
$$\lbrace \hat{u},\hat{v}\rbrace:=\pi_{n+m-1}(uv-\chi(u,v)vu)$$
with $u\in \mathcal{U}(\mathfrak{g})_n$ and $v\in \mathcal{U}(\mathfrak{g})_m$ and $\hat{u}=\pi_n(u)$, $\hat{v}=\pi_m(v)$
\end{thm}
\begin{proof} The bracket $\lbrace \hat{u},\hat{v}\rbrace$ is independent of the choose of $u$ and $v$.  If $\hat{X},\hat{Y}\in gr^1(\mathcal{U}(\mathfrak{g}))$, then the bracket satisfies the properties of Poisson Color. This bracket extends on $gr(\mathcal{U}(\mathfrak{g}))$ by belinearity. The only property we need to verify is the Leibniz property . Let $\hat{X},\hat{Y},\hat{Z}\in gr^1(\mathcal{U}(\mathfrak{g}))$ and $X,Y,Z\in\mathcal{U}(\mathfrak{g})_1$  we have  

$$\lbrace \hat{X},\hat{Y}\hat{Z}\rbrace= \pi_{1+2-1}(X(YZ)-\omega(X,YZ)(YZ)X)$$

$$=\pi_2((XY-\omega(X,Y)YX)Z+\omega(X,Y)Y(XZ-\omega(X,Z)ZX))$$
$$=((\pi_1(X)\pi_1(Y)-\omega(X,Y)\pi_1(Y)\pi(X)\pi_1(Z)+\omega(X,Y)\pi_1(Y)(\pi_1(X)\pi(Z)-\omega(X,Z)\pi_1(Z)\pi_1(X))$$
$$=\lbrace \hat{X},\hat{Y}\rbrace\hat{Z}+\omega(\hat{X},\hat{Y})\lbrace \hat{X},\hat{Z}\rbrace$$
Since $\mathcal{U}(\mathfrak{g})_1=\korps\oplus\mathfrak{g}$, we have $\omega(\hat{X},\hat{Y})=\omega(X,Y)$.\\
 By Poincar\'{e}-Birkhoff-Witt Theorem, this bracket regives the Color Lie algebra of $\mathfrak{g}$.
\end{proof}

\begin{de}Let $(A,m_0,\lbrace,\rbrace)$ be a Color Poisson bracket. A graded star product of $A$ is graded associative algebra deformation $(A[[t]],m_t)$ of $A$ such that 
$$m_1(a_{\alpha},a_{\beta})-\chi(\alpha,\beta)m_1(a_{\beta}, a_{\alpha})=\lbrace a_{\alpha},a_{\beta}\rbrace $$
for $a_{\alpha}\in A_{\alpha}, a_{\beta}\in A_{\beta}$
\end{de}
\begin{thm}Poincar\'{e}-Birkhoff-Witt. Let $(\mathfrak{g},[,])$ be a Color Lie algebra, $\mathcal{U}(\mathfrak{g})$ its universal enveloping algebra and the associated graded algebra $gr(\mathcal{U}(\mathfrak{g}))$. The symmetrization map $\omega: gr(\mathcal{U}(\mathfrak{g}))\to \mathcal{U}(\mathfrak{g})$  is homogenous isomorphism of degree zero of $\mathfrak{g}$ module, defined by
$$ \omega_p(X_{1}\bullet\cdots\bullet X_{p}):=\sum_{\sigma\in \mathcal{S}_p}\prod _{(i,j):\sigma(j)<\sigma(i)}\chi(i,j)X_{\sigma(1)}\circ\cdots\circ X_{\sigma(p)}$$
where $\mathcal{S}_p$ is the group of permutation of ordr $p$ and the product is extended over all $r,s\in {1,2,...,n}$ such that $r<s$ and $\pi^{-1}(r)>\pi^{-1}(s)$, $\bullet$ is the multiplication of $gr(\mathcal{U}(\mathfrak{g}))$ and $\circ$ is the multiplication of $\mathcal{U}(\mathfrak{g})$. We have 
$$\mathcal{U}(\mathfrak{g})=\oplus_{j\geqslant 0}\omega(gr^j(\mathcal{U}(\mathfrak{g}))$$
\end{thm}
\begin{thm}Let $(\mathfrak{g},[,])$ be a Color Lie algebra, $\mathcal{U}(\mathfrak{g})$ its universal enveloping algebra and the associated graded algebra $gr(\mathcal{U}(\mathfrak{g}))$. The symmetrization map $\omega$ defines a star product over the Color Poisson algebra $gr(\mathcal{U}(\mathfrak{g}))$ given by
$$u\star_t v=\sum_{n\geqslant 0}t^n\omega^{-1}((\omega(u)\circ \omega(v))_{p+q-n})=\sum^{p+q-1}_{n= 0}t^n\omega^{-1}((\omega(u)\circ \omega(v))_{p+q-n})$$
where $u\in gr^p(\mathcal{U}(\mathfrak{g})$ and $v\in gr^q(\mathcal{U}(\mathfrak{g})$ .
\end{thm}

\end{document}